\documentclass[oneside,english]{amsart}
\usepackage[T1]{fontenc}
\usepackage[latin9]{inputenc}
\usepackage{url}
\usepackage{amssymb}
\usepackage{esint}

\makeatletter
\numberwithin{equation}{section} 
\numberwithin{figure}{section} 
  \theoremstyle{plain}
  \newtheorem{thm}{Theorem}[section]
  \theoremstyle{plain}
  \newtheorem{cor}[thm]{Corollary}
  \theoremstyle{plain}
  \newtheorem{lem}[thm]{Lemma}
  \theoremstyle{remark}
  \newtheorem{rem}[thm]{Remark}
  \theoremstyle{remark}
  \newtheorem*{acknowledgement*}{Acknowledgement}

\usepackage{babel}
\makeatother

\begin{document}

\title{Explicit Inversions of Certain Matrices II}

\author{Ruiming Zhang}

\begin{abstract}
In this note, we apply kernel polynomials to find the explicit inverses
for some some Hankel matrices associated with $q$-orthogonal polynomials.
\end{abstract}


\curraddr{School of Mathematical Sciences\\
Guangxi Normal University\\
Guilin City, Guangxi 541004\\
P. R. China.}


\email{ruimingzhang@yahoo.com}

\maketitle

\section{Introduction}

In the theory of orthogonal polynomials, We could calculate the determinants
of some Hankel matrices once we know the three term recurrence relation
for the associated orthogonal polynomials and vice versa. It is well-known
that the kernel polynomials of the orthogonal polynomials encode important
information about the associated Hankel matrices. These matrices are
generalizations of the Hilbert matrices. In this note we invert some
Hankel matrices associated with some $q$-polynomials by using the
kernel polynomials. 

The rest of the introduction could be found in \cite{Zhang}, we copy
it here for the convenience of the reader.

\begin{thm}
\label{thm:1.1}Given a probability measure $\mu$ on $\mathbb{R}$
with a support of infinite many points. Let us consider the Hilbert
space of $\mu$-measurable functions \begin{equation}
\mathcal{X}:=\left\{ f(x)|\int|f(x)|^{2}d\mu(x)<\infty\right\} \label{eq:1.1}\end{equation}
with the inner product defined as\begin{equation}
(f,g):=\int f(x)\overline{g(x)}d\mu(x),\quad f,g\in\mathcal{X}.\label{eq:1.2}\end{equation}
Assume that $\left\{ w_{n}(x)\right\} _{n=0}^{\infty}$ is a sequence
of linearly independent functions in $\mathcal{X}$ with $w_{0}(x)=1$.
Let

\begin{equation}
\alpha_{jk}:=\int w_{j}(x)\overline{w_{k}(x)}d\mu(x),\quad j,k=0,1...,\label{eq:1.3}\end{equation}
 \begin{equation}
\Pi_{n}:=\left(\begin{array}{cccc}
\alpha_{00} & \alpha_{01} & \dots & \alpha_{0n}\\
\alpha_{10} & \alpha_{11} & \dots & \alpha_{1n}\\
\vdots & \vdots & \vdots & \vdots\\
\alpha_{n0} & \alpha_{n1} & \dots & \alpha_{nn}\end{array}\right),\quad n\in\mathbb{N}\cup\left\{ 0\right\} ,\label{eq:1.4}\end{equation}
 and

\begin{equation}
\Delta_{n}:=\det\Pi_{n},\quad n\in\mathbb{N}\cup\left\{ 0\right\} .\label{eq:1.5}\end{equation}
Then the $n$-th orthonormal function with positive coefficient in
$w_{n}(x)$ is given by the formula

\begin{equation}
p_{n}(x)=\frac{1}{\sqrt{\Delta_{n}\Delta_{n-1}}}\det\left(\begin{array}{ccccc}
\alpha_{00} & \alpha_{01} & \alpha_{02} & \dots & \alpha_{0n}\\
\alpha_{10} & \alpha_{11} & \alpha_{12} & \dots & \alpha_{1n}\\
\vdots & \vdots & \vdots & \ddots & \vdots\\
\alpha_{n0} & \alpha_{n1} & \alpha_{n2} & \dots & \alpha_{nn}\\
w_{0}(x) & w_{1}(x) & w_{2}(x) & \dots & w_{n}(x)\end{array}\right)\label{eq:1.6}\end{equation}
 for $n\in\mathbb{N}$ , with\begin{equation}
p_{0}(x)=w_{0}(x)=1.\label{eq:1.7}\end{equation}
 Furthermore, the coefficient of $p_{n}(x)$ in $w_{n}(x)$ is \begin{equation}
\gamma_{n}:=\sqrt{\frac{\Delta_{n-1}}{\Delta_{n}}}.\label{eq:1.8}\end{equation}

\end{thm}
\begin{proof}
The proof is the same as for the case $w_{n}(x)=x^{n}$, which could
be found in any orthogonal polynomials textbooks such as \cite{Szego}.
\end{proof}
\begin{cor}
\label{cor:1.1} For $n\in\mathbb{N}$, we have \begin{equation}
\Delta_{n}=\prod_{k=1}^{n}\frac{1}{\gamma_{n}^{2}}.\label{eq:1.9}\end{equation}

\end{cor}
\begin{proof}
This is a trivial consequence of \eqref{eq:1.7} and \eqref{eq:1.8}. 
\end{proof}
\begin{lem}
\label{lem:1.1} Let \begin{equation}
k_{n}(x,y):=\sum_{k=0}^{n}p_{k}(x)\overline{p_{k}(y)},\quad n\in\mathbb{N}\cup\left\{ 0\right\} .\label{eq:1.10}\end{equation}
 Then, for any $\pi(x)$ in the linear span of $\left\{ w_{k}(x)\right\} _{0}^{n}$
, we have

\begin{equation}
\int\pi(x)\overline{k_{n}(x,y)}d\mu(x)=\pi(y).\label{eq:1.11}\end{equation}
 
\end{lem}
\begin{proof}
To see \eqref{eq:1.11}, just expand $\pi(x)$ in $p_{k}(x),k=0,1\dots,n$. 
\end{proof}
\begin{lem}
\label{lem:1.2} For each $n\in\mathbb{N}\cup\left\{ 0\right\} $,
the function $k_{n}(x,y)$ satisfying \eqref{eq:1.11} is unique.
\end{lem}
\begin{proof}
Suppose there are two such functions $h_{n}(x,y)$ and $k_{n}(x,y)$
, then, 

\begin{align}
0 & <||h_{n}(\cdot,y)-k_{n}(\cdot,y)||^{2}\label{eq:1.12}\\
= & (h_{n}(\cdot,y)-k_{n}(\cdot,y),h_{n}(\cdot,y)-k_{n}(\cdot,y))\nonumber \\
= & (h_{n}(\cdot,y)-k_{n}(\cdot,y),h_{n}(\cdot,y))-(h_{n}(\cdot,y)-k_{n}(\cdot,y),k_{n}(\cdot,y))\nonumber \\
= & 0,\nonumber \end{align}
 which is a contradiction. 
\end{proof}
\begin{lem}
\label{lem:1.3} Let \begin{align}
(\beta_{jk})_{0\le j,k\le n}: & =\Pi_{n}^{-1},\quad n\in\mathbb{N}\cup\left\{ 0\right\} .\label{eq:1.13}\end{align}
Then,\begin{align}
k_{n}(x,y) & =\sum_{j,k=0}^{n}\beta_{jk}\overline{w_{j}(y)}w_{k}(x).\label{eq:1.14}\end{align}

\end{lem}
\begin{proof}
Let \begin{align}
f(x) & =\sum_{k=0}^{n}u_{k}w_{k}(x),\label{eq:1.15}\end{align}
 then,

\begin{align}
 & (f(\cdot),\sum_{j,k=0}^{n}\beta_{jk}\overline{w_{j}(y)}w_{k}(\cdot))\label{eq:1.16}\\
= & \sum_{m=0}^{n}u_{m}(w_{m}(\cdot),k(\cdot,y))\nonumber \\
= & \sum_{m=0}^{n}u_{m}\sum_{j,k=0}^{n}\overline{\beta_{jk}}w_{j}(y)(w_{m},w_{k})\nonumber \\
= & \sum_{m=0}^{n}u_{m}\sum_{j=0}^{n}w_{j}(y)\sum_{k=0}^{n}\overline{\beta_{jk}\alpha_{km}}\nonumber \\
= & f(y).\nonumber \end{align}
 By Lemma \ref{lem:1.2}, we have\begin{align}
k_{n}(x,y) & =\sum_{j,k=0}^{n}\beta_{jk}\overline{w_{j}(y)}w_{k}(x).\label{eq:1.17}\end{align}

\end{proof}
\begin{cor}
\label{cor:1.3}The kernel in Lemma \ref{lem:1.1} is also given by

\begin{align}
k_{n}(x,y) & =-\frac{1}{\Delta_{n}}\det\left(\begin{array}{ccccc}
0 & 1 & \overline{w_{1}(y)} & \cdots & \overline{w_{n}(y)}\\
1 & \alpha_{00} & \alpha_{01} & \dots & \alpha_{0n}\\
w_{1}(x) & \alpha_{10} & \alpha_{11} & \dots & \alpha_{1n}\\
\vdots & \vdots & \vdots & \ddots & \vdots\\
w_{n}(x) & \alpha_{n0} & \alpha_{n1} & \dots & \alpha_{nn}\end{array}\right)\label{eq:1.18}\end{align}
 for $n\in\mathbb{N}\cup\left\{ 0\right\} $. 
\end{cor}
\begin{proof}
Since \begin{align}
k_{n}(x,y) & =\sum_{j,k=0}^{n}\beta_{jk}\overline{w_{j}(y)}w_{k}(x),\label{eq:1.19}\end{align}
 with \begin{align}
\left(\beta_{jk}\right)_{0\le i,j\le n} & =\Pi_{n}^{-1}.\label{eq:1.20}\end{align}
 Then,\begin{align}
\beta_{jk} & =\frac{\Pi_{n}(k,j)}{\det\Pi_{n}}=\frac{\Pi_{n}(k,j)}{\Delta_{n}},\label{eq:1.21}\end{align}
 where $\Pi_{n}(k,j)$ is the $(k,j)$-th co-factor. Therefore,\begin{align}
k_{n}(x,y) & =\frac{1}{\Delta_{n}}\sum_{j,k=0}^{n}\Pi_{n}(k,j)\overline{w_{j}(y)}w_{k}(x).\label{eq:1.22}\end{align}
 It is clear that\begin{align}
\sum_{j,k=0}^{n}\Pi_{n}(k,j)\overline{w_{j}(y)}w_{k}(x) & =-\left|\begin{array}{cc}
0 & \mathbf{(\overline{W(y)})}^{T}\\
\mathbf{W(x)} & \Pi_{n}\end{array}\right|,\label{eq:1.23}\end{align}
 by direct determination expansion, which is \begin{align}
k_{n}(x,y) & =-\frac{1}{\Delta_{n}}\left|\begin{array}{cc}
0 & \mathbf{(\overline{W(y)})}^{T}\\
\mathbf{W(x)} & \Pi_{n}\end{array}\right|,\label{eq:1.24}\end{align}
where\begin{align}
\mathbf{W(x)} & =\left(\begin{array}{c}
1\\
w_{1}(x)\\
\vdots\\
w_{n}(x)\end{array}\right),\label{eq:1.25}\end{align}
and\begin{align}
\mathbf{\mathbf{(\overline{W(y)})}^{T}} & =\left(\begin{array}{cccc}
1, & \overline{w_{1}(y)} & ,\cdots, & \overline{w_{n}(y)}\end{array}\right).\label{eq:1.26}\end{align}

\end{proof}
Lemma \ref{lem:1.3} enables us to compute the inverse the Gram matrix
in terms of the orthonormal functions $\left\{ p_{n}(x)\right\} _{n=0}^{\infty}$.

\begin{cor}
\label{cor:1.4} Assume that $\left\{ w_{n}(x)\right\} _{n=0}^{\infty}$,
$\left\{ p_{n}(x)\right\} _{n=0}^{\infty}$ and $\Pi_{n}=(\alpha_{jk})_{0\le j,k\le n}$
as in Theorem \ref{thm:1.1}. Suppose we have two families of linear
functionals $\left\{ u_{k}\right\} _{k=0}^{\infty}$ and $\left\{ v_{k}\right\} _{k=0}^{\infty}$
over the linear space generated by $\left\{ w_{n}(x)\right\} _{n=0}^{\infty}$
with\begin{align}
u_{j}(w_{k}) & =\delta_{jk},\label{eq:1.27}\end{align}
 and

\begin{align}
v_{j}(\overline{w_{k}}) & =\delta_{jk}\label{eq:1.28}\end{align}
 for $j,k=0,1,...$. Then, 

\begin{align}
\beta_{jk} & =\sum_{m=0}^{n}u_{k}(p_{m}(x))v_{j}(\overline{p_{m}(y)}),\label{eq:1.29}\end{align}
where \begin{equation}
(\beta_{jk})_{0\le j,k\le n}=\Pi_{n}^{-1}.\label{eq:1.30}\end{equation}

\end{cor}
\begin{proof}
From Lemma \ref{lem:1.3}, we have\begin{equation}
\sum_{j,k=0}^{n}\beta_{jk}\overline{w_{j}(y)}w_{k}(x)=\sum_{m=0}^{n}\overline{p_{m}(y)}p_{m}(x).\label{eq:1.31}\end{equation}
Then we apply the functional $u_{j}$ and $v_{k}$ both sides of the
above equation, the claim of the corollary follows. 
\end{proof}

\section{Main Results\label{sec:Main-Results}}

As usual, \cite{Andrews,Koekoek},

\begin{align}
(a;q)_{\infty}: & =\prod_{m=0}^{\infty}(1-aq^{m}),\quad a\in\mathbb{C},\quad q\in(0,1),\label{eq:2.1}\end{align}
and\begin{align}
(a;q)_{m}: & =\frac{(a;q)_{\infty}}{(aq^{m};q)_{\infty}},\quad m\in\mathbb{Z}.\label{eq:2.2}\end{align}
 We will use the following short hand notations\begin{align}
\left[\begin{array}{c}
m\\
j\end{array}\right]_{q}: & =\frac{(q;q)_{m}}{(q;q)_{j}(q;q)_{m-j}},\quad j\le m,\quad j,m\in\mathbb{N}\cup\left\{ 0\right\} ,\label{eq:2.3}\end{align}
 and

\begin{align}
(a_{1},a_{2},...,a_{n};q)_{m}: & =\prod_{k=1}^{n}(a_{k};q)_{m},\quad n\in\mathbb{N},\quad m\in\mathbb{Z},\quad a_{1},\dots,a_{n}\in\mathbb{C}.\label{eq:2.4}\end{align}
The basic hypergeometric function ${}_{r}\phi_{s}$ with complex parameters
$a_{1},...,a_{r};b_{1},...,b_{s}$ is formally defined as\begin{align}
{}_{r}\phi_{s}\left(\begin{array}{c}
a_{1},...,a_{r}\\
b_{1},...,b_{s}\end{array};q;z\right) & :=\sum_{m=0}^{\infty}\frac{(a_{1},...,a_{r};q)_{m}}{(q,b_{1},...,b_{s};q)_{m}}\left((-1)^{n}q^{(n-1)n/2}\right)^{s+1-r}z^{m}.\label{eq:2.5}\end{align}
The $q$-binomial theorem is the following summation formula\begin{align}
\frac{(az;q)_{\infty}}{(z;q)_{\infty}} & =\sum_{n=0}^{\infty}\frac{(a;q)_{n}}{(q;q)_{n}}z^{n}\quad|z|<1,\label{eq:2.6}\end{align}
and the Ramanujan ${}_{1}\psi_{1}$ summation formula is \begin{align}
\sum_{n=-\infty}^{\infty}\frac{(a;q)_{n}}{(b;q)_{n}}z^{n} & =\frac{(b/a,q,q/az,az;q)_{\infty}}{(b,b/az,q/a,z;q)_{\infty}},\quad|\frac{b}{a}|<|z|<1.\label{eq:2.7}\end{align}

\begin{thm}
\label{thm:2.7}For $n\in\mathbb{N}\cup\left\{ 0\right\} $, the matrix
\begin{align}
 & \left(\frac{(aq;q)_{j+k}}{(abq^{2};q)_{j+k}}\right)_{j,k=0}^{n}\label{eq:2.8}\end{align}
 has determinant
\end{thm}
\begin{align}
 & \det\left(\frac{(aq;q)_{j+k}}{(abq^{2};q)_{j+k}}\right)_{j,k=0}^{n}\label{eq:2.9}\\
= & a^{n(n+1)/2}q^{n(n+1)(2n+1)/6}\prod_{k=1}^{n}\frac{(q,aq,bq,abq;q)_{k}}{(abq,abq^{2},abq^{2},abq^{3};q^{2})_{k}}.\nonumber \end{align}
 The matrix \eqref{eq:2.8} is invertible for \begin{align}
q^{k},aq^{k},bq^{k},abq^{k}\neq1,\quad k & \in\mathbb{N},\quad a,q\neq0,\label{eq:2.10}\end{align}
 its inverse matrix $(\gamma_{jk})_{j,k=0}^{n}$ has elements\begin{align}
\gamma_{jk} & =\frac{(-1)^{j+k}q^{[j(j+1)+k(k+1)]/2}}{(aq;q)_{j}(aq;q)_{k}}\label{eq:2.11}\\
\times & \sum_{m=\max(j,k)}^{n}\frac{(1-abq^{2m+1})(aq,abq;q)_{m}}{(1-abq)(q,bq;q)_{m}(aq^{j+k+1})^{m}}\nonumber \\
\times & \left[\begin{array}{c}
m\\
j\end{array}\right]_{q}\cdot\left[\begin{array}{c}
m\\
k\end{array}\right]_{q}(abq^{m+1};q)_{j}(abq^{m+1};q)_{k}.\nonumber \end{align}

\begin{rem}
Let $a=b=1$ in \eqref{eq:2.8}, then take $q\to1$ we would get the
classical Hilbert matrix. 
\end{rem}
\begin{thm}
\label{thm:2.9}For $n\in\mathbb{N}\cup\left\{ 0\right\} $ and $q\neq0$,
the matrix\begin{align}
 & \left((q^{\alpha+1};q)_{j+k}q^{-\binom{j+k}{2}}\right)_{j,k=0}^{n}\label{eq:2.12}\end{align}
 has determinant

\begin{align}
\det\left((q^{\alpha+1};q)_{j+k}q^{-\binom{j+k}{2}}\right)_{j,k=0}^{n} & =\frac{\prod_{k=1}^{n}(q,q^{\alpha+1};q)_{k}}{q^{n(n+1)(4n-1)/6}}.\label{eq:2.13}\end{align}
 The matrix \eqref{eq:2.12} is invertible for \begin{align}
-\alpha\notin\mathbb{N},\quad q^{k}\neq1,k & \in\mathbb{N},\label{eq:2.14}\end{align}
 the inverse matrix $(\gamma_{jk})_{j,k=0}^{n}$ has elements \begin{align}
\gamma_{jk} & =\frac{q^{[j(j-1)+k(k-1)]}}{(q^{\alpha+1};q)_{j}(q^{\alpha+1};q)_{k}}\sum_{m=\max(j,k)}^{n}\frac{(q^{\alpha+1};q)_{m}}{(q;q)_{m}}\left[\begin{array}{c}
m\\
j\end{array}\right]_{q}\left[\begin{array}{c}
m\\
k\end{array}\right]_{q}q^{m}.\label{eq:2.15}\end{align}
 
\end{thm}
In our last result, we use the following notations \begin{align}
x & =\left\lfloor x\right\rfloor +\left\{ x\right\} ,\quad x\in\mathbb{R},\quad\left\lfloor x\right\rfloor \in\mathbb{Z},\quad\left\{ x\right\} \in[0,1).\label{eq:2.16}\end{align}

\begin{thm}
\label{thm:2.10}For $n\in\mathbb{N}\cup\left\{ 0\right\} $ and $q\neq0$,
the matrix\begin{align}
 & \left(\frac{\left[1+(-1)^{j+k}\right]}{2}q^{-(j+k)^{2}/4}(q;q^{2})_{\left\lfloor \frac{j+k}{2}\right\rfloor }\right)_{j,k=0}^{n}\label{eq:2.17}\end{align}
has the determinant\begin{align}
\det\left(\frac{\left[1+(-1)^{j+k}\right]}{2}q^{-(j+k)^{2}/4}(q;q^{2})_{\left\lfloor \frac{j+k}{2}\right\rfloor }\right)_{j,k=0}^{n} & =\frac{\prod_{k=0}^{n}(q;q)_{k}}{q^{n(n+1)(2n+1)/6}}.\label{eq:2.18}\end{align}
For \begin{align}
q^{k}\neq1,\, k & \in\mathbb{N},\label{eq:2.19}\end{align}
the matrix \eqref{eq:2.17} is invertible and has inverse \begin{align}
 & \left(\sum_{m=\max(j,k)}^{n}\frac{\left[\begin{array}{c}
m\\
j\end{array}\right]_{q}\left[\begin{array}{c}
m\\
k\end{array}\right]_{q}(q,q^{2})_{\left\lfloor \frac{m-j}{2}\right\rfloor }(q,q^{2})_{\left\lfloor \frac{m-k}{2}\right\rfloor }\cos\frac{(m-j)\pi}{2}\cos\frac{(m-k)\pi}{2}}{(-1)^{m}(q;q)_{m}q^{-m-\binom{j}{2}-\binom{k}{2}}}\right)_{j,k=0}^{n}.\label{eq:2.20}\end{align}

\end{thm}

\section{Proofs }

Given a polynomial $f(x)$, we define the $q$-difference operator
as \cite{Andrews,Koekoek}\begin{align}
(D_{q}f)(x): & =\frac{f(x)-f(qx)}{(1-q)x}.\label{eq:3.1}\end{align}
 Let $w_{k}(x):=x^{k}$ for $k=0,1,...$. then we have\begin{align}
D_{q}w_{n}(x) & =\frac{1-q^{n}}{1-q}w_{n-1}(x),\label{eq:3.2}\end{align}
 and

\begin{align}
\left[D_{q}^{k}w_{n}(x)\right]_{x=0} & =\frac{(q;q)_{n}}{(1-q)^{n}}\delta_{kn}.\label{eq:3.3}\end{align}
 For this polynomial sequence we take \begin{align}
u_{k}(f(x)): & =v_{k}(f(x)):=\frac{(1-q)^{k}}{(q;q)_{k}}\left[(D_{q}^{k}f)(x)\right]_{x=0},\label{eq:3.4}\end{align}
 where $f(x)$ is a polynomial in variable $x$.

\subsection{Proof for Theorem \ref{thm:2.7} }

The little $q$-Jacobi polynomials $\left\{ p_{n}(x;a,b|q)\right\} _{n=0}^{\infty}$
are defined as \cite{Andrews,Koekoek} \begin{align}
p_{n}(x;a,b|q): & ={}_{2}\phi_{1}\left(\begin{array}{c}
q^{-n},abq^{n+1}\\
aq\end{array};q;qx\right),\quad0<q,aq,bq<1\label{eq:3.5}\end{align}
 for $n\ge0$, and we assume that\begin{align}
p_{-1}(x;a,b|q): & =0,\quad p_{0}(x;a,b|q):=1.\label{eq:3.6}\end{align}
 The little $q$-Jacobi polynomials $\left\{ p_{n}(x;a,b|q)\right\} _{n=0}^{\infty}$
have the orthogonal relation

\begin{align}
\sum_{k=0}^{\infty}\frac{(bq;q)_{k}(aq)^{k}}{(q;q)_{k}}p_{m}(q^{k};a,b|q)p_{n}(q^{k};a,b|q) & =h_{n}(a,b|q)\delta_{mn}\label{eq:3.7}\end{align}
 for $m,n\ge0$ with\begin{align}
h_{n}(a,b|q): & =\frac{(abq^{2};q)_{\infty}}{(aq;q)_{\infty}}\frac{(1-abq)(aq)^{n}}{(1-abq^{2n+1})}\frac{(q,bq;q)_{n}}{(aq,abq;q)_{n}}.\label{eq:3.8}\end{align}
 The little $q$-Jacobi polynomials satisfy the following difference
relation\begin{align}
D_{q}^{k}p_{n}(x;a,b|q) & =\frac{q^{k}(q^{-n};q)_{k}(abq^{n+1};q)_{k}}{(1-q)^{k}(aq;q)_{k}}p_{n-k}(x;aq^{k},bq^{k}|q),\label{eq:3.9}\end{align}
 for $n,k=0,1...$ and $n\ge k$. The orthonormal polynomial\begin{align}
p_{n}(x) & =\frac{(-1)^{n}p_{n}(x;aq,bq|q)}{\sqrt{h_{n}(a,b|q)}}\label{eq:3.10}\end{align}
 have the leading coefficient\begin{align}
\gamma_{n} & =\sqrt{\frac{(aq;q)_{\infty}(abq,abq^{2},abq^{2},abq^{3};q^{2})_{n}}{(abq^{2};q)_{\infty}(q,aq,bq,abq;q)_{n}a^{n}q^{n^{2}}}}\label{eq:3.11}\end{align}
 for $n=0,1...$.

The moments are given by the formula\begin{align}
\mu_{n} & =\sum_{m=0}^{\infty}\frac{(bq;q)_{m}(aq)^{m}q^{nm}}{(q;q)_{m}},\label{eq:3.12}\end{align}
 or\begin{equation}
\mu_{n}=\frac{(abq^{n+2};q)_{\infty}}{(aq^{n+1};q)_{\infty}}.\label{eq:3.13}\end{equation}
 by using the $q$-binomial theorem. 

For any $n=0,1,...$, the matrix 

\begin{align}
H_{n} & =\left(\frac{(abq^{j+k+2};q)_{\infty}}{(aq^{j+k+1};q)_{\infty}}\right)_{j,k=0}^{n}\label{eq:3.14}\end{align}
 has the determinant\begin{align}
\det\left(\frac{(abq^{j+k+2};q)_{\infty}}{(aq^{j+k+1};q)_{\infty}}\right)_{j,k=0}^{n} & =\left(\frac{(abq^{2};q)_{\infty}a^{n/2}q^{n(2n+1)/6}}{(aq;q)_{\infty}}\right)^{n+1}\label{eq:3.15}\\
\times & \prod\frac{(q,aq,bq,abq;q)_{k}}{(abq,abq^{2},abq^{2},abq^{3};q^{2})_{k}},\nonumber \end{align}
 or\begin{align}
\det\left(\frac{(aq;q)_{j+k}}{(abq^{2};q)_{j+k}}\right)_{j,k=0}^{n} & =a^{n(n+1)/2}q^{n(n+1)(2n+1)/6}\label{eq:3.16}\\
\times & \prod_{k=1}^{n}\frac{(q,aq,bq,abq;q)_{k}}{(abq,abq^{2},abq^{2},abq^{3};q^{2})_{k}}\nonumber \end{align}
by applying Theorem \ref{thm:1.1}.

Under the condition $0<q,aq,bq<1$, the elements of $H_{n}^{-1}=(\beta_{jk})_{0\le j,k\le n}$
are given by

\begin{align}
\beta_{jk} & =\frac{q^{j+k}(aq;q)_{\infty}}{(abq^{2};q)_{\infty}}\sum_{m=\max(j,k)}^{n}\frac{(1-abq^{2m+1})(aq,abq;q)_{m}}{(1-abq)(aq)^{m}(q,bq;q)_{m}}\label{eq:3.17}\\
\times & \frac{(q^{-m};q)_{j}(abq^{m+1};q)_{j}}{(q;q)_{j}(aq;q)_{j}}\frac{(q^{-m};q)_{k}(abq^{m+1};q)_{k}}{(q;q)_{j}(aq;q)_{k}}\nonumber \end{align}
 for $j,k=0,1,...,n$. Then, the matrix \begin{align}
 & \left(\frac{(aq;q)_{j+k}}{(abq^{2};q)_{j+k}}\right)_{j,k=0}^{n}\label{eq:3.18}\end{align}
 has the inverse matrix $(\gamma_{jk})_{j,k=0}^{n}$ with\begin{align}
\gamma_{jk} & =\frac{(-1)^{j+k}q^{[j(j+1)+k(k+1)]/2}}{(aq;q)_{j}(aq;q)_{k}}\label{eq:3.19}\\
\times & \sum_{m=\max(j,k)}^{n}\frac{(1-abq^{2m+1})(aq,abq;q)_{m}}{(1-abq)(q,bq;q)_{m}(aq^{j+k+1})^{m}}\nonumber \\
\times & \left[\begin{array}{c}
m\\
j\end{array}\right]_{q}\left[\begin{array}{c}
m\\
k\end{array}\right]_{q}(abq^{m+1};q)_{j}(abq^{m+1};q)_{k}.\nonumber \end{align}
Since the formulas \eqref{eq:3.16}, \eqref{eq:3.18} and \eqref{eq:3.19}
contain only rational functions of $a$, $b$ and $q^{k},k\in\mathbb{N}$,
the original restrictions in \eqref{eq:3.5} could be dropped and
Theorem \ref{thm:2.7} follows.

\subsection{Proof for Theorem \ref{thm:2.9} }

The $q$-Laguerre polynomials $\left\{ L_{n}^{(\alpha)}(x;q)\right\} _{n=0}^{\infty}$
are defined as \cite{Andrews,Koekoek}\begin{align}
L_{n}^{(\alpha)}(x;q): & =\frac{(q^{\alpha+1};q)_{n}}{(q;q)_{n}}\sum_{k=0}^{n}\frac{(q^{-n};q)_{k}q^{\binom{k+1}{2}}(-x)^{k}q^{(\alpha+n)k}}{(q;q)_{k}(q^{\alpha+1};q)_{k}},\quad q\in(0,1),\,\alpha>-1\label{eq:3.20}\end{align}
 for $n\ge0$, and we assume that\begin{align}
L_{-1}^{(\alpha)}(x;q): & =0,\quad L_{0}^{(\alpha)}(x;q):=1.\label{eq:3.21}\end{align}
The $q$-Laguerre polynomials have an orthogonal relation\begin{align}
\sum_{k=-\infty}^{\infty}\frac{q^{k\alpha+k}}{(-cq^{k};q)_{\infty}}L_{m}^{(\alpha)}(cq^{k};q)L_{n}^{(\alpha)}(cq^{k};q) & =h_{n}(c,\alpha|q)\delta_{mn}\label{eq:3.22}\end{align}
 for $m,n=0,1,...$ with\begin{align}
h_{n}(c,\alpha|q): & =\frac{(q,-cq^{\alpha+1},-c^{-1}q^{-\alpha};q)_{\infty}}{(q^{\alpha+1},-c,-c^{-1}q;q)_{\infty}}\frac{(q^{\alpha+1};q)_{n}}{(q;q)_{n}q^{n}}.\label{eq:3.23}\end{align}
 The $q$-Laguerre polynomials satisfy the following difference relation\begin{align}
D_{q}^{k}L_{n}^{(\alpha)}(x;q) & =\frac{q^{k\alpha+k^{2}}}{(1-q)^{k}}L_{n-k}^{(\alpha+k)}(q^{k}x;q)\label{eq:3.24}\end{align}
 for $n,k=0,1...$. The orthonormal $q$-Laguerre polynomial\begin{align}
\ell_{n}(x): & =\frac{L_{n}^{(\alpha)}(x;q)}{\sqrt{h_{n}(c,\alpha|q)}}\label{eq:3.25}\end{align}
 has the leading coefficient\begin{align}
\gamma_{n} & =\frac{q^{(\alpha+n+1/2)n}}{\sqrt{(q,q^{\alpha+1};q)_{n}}}\sqrt{\frac{(q^{\alpha+1},-c,-c^{-1}q;q)_{\infty}}{(q,-cq^{\alpha+1},-c^{-1}q^{-\alpha};q)_{\infty}}}\label{eq:3.26}\end{align}
 for $n=0,1...$.

The moments are given by \begin{align}
\mu_{n} & =\sum_{k=-\infty}^{\infty}\frac{c^{n}q^{k(\alpha+n+1)}}{(-cq^{k};q)_{\infty}}\label{eq:3.27}\\
= & \frac{c^{n}}{(-c;q)_{\infty}}\sum_{k=-\infty}^{\infty}(-c;q)_{k}q^{k(\alpha+n+1)},\nonumber \end{align}
 which could be summed by the formula \eqref{eq:2.7},\begin{align}
\mu_{n} & =\frac{c^{n}}{(-c;q)_{\infty}}\frac{(q,-cq^{\alpha+n+1},-q^{-\alpha-n}/c;q)_{\infty}}{(-q/c,q^{\alpha+n+1};q)_{\infty}}\label{eq:3.28}\\
= & \frac{c^{n}(q^{\alpha+1};q)_{n}(-q^{-\alpha-n}/c;q)_{n}}{(-cq^{\alpha+1};q)_{n}}\frac{(q,-cq^{\alpha+1},-q^{-\alpha}/c;q)}{(-c,-q/c,q^{\alpha+1};q)_{\infty}}\nonumber \\
= & \frac{(q^{\alpha+1};q)_{n}}{q^{n\alpha+n(n+1)/2}}\frac{(q,-cq^{\alpha+1},-q^{-\alpha}/c;q)}{(-c,-q/c,q^{\alpha+1};q)_{\infty}}.\nonumber \end{align}
 Then for any $n=0,1,...$, the matrix

\begin{align}
H_{n} & =\left(\frac{(q^{\alpha+1};q)_{j+k}}{q^{\left[(j+k)\alpha+(j+k)(j+k+1)/2\right]}}\frac{(q,-cq^{\alpha+1},-q^{-\alpha}/c;q)}{(-c,-q/c,q^{\alpha+1};q)_{\infty}}\right)_{j,k=0}^{n}\label{eq:3.29}\end{align}
 has the determinant\begin{align}
\det H_{n} & =\left(\frac{(q,-cq^{\alpha+1},-c^{-1}q^{-\alpha};q)_{\infty}}{(q^{\alpha+1},-c,-c^{-1}q;q)_{\infty}}\right)^{n+1}\prod_{k=0}^{n}\frac{(q,q^{\alpha+1};q)_{k}}{q^{k\alpha+k(k+1)/2}},\label{eq:3.30}\end{align}
 or\begin{align}
\det\left((q^{\alpha+1};q)_{j+k}q^{-\binom{j+k}{2}}\right)_{j,k=0}^{n} & =\frac{\prod_{k=1}^{n}(q,q^{\alpha+1};q)_{k}}{q^{n(n+1)(4n-1)/6}},\label{eq:3.31}\end{align}
by applying Theorem \ref{thm:1.1}. By applying Corollary \ref{cor:1.6},
we obtain the inverse matrix $H_{n}^{-1}=(\beta_{jk})_{j,k=0}^{n}$
with\begin{align}
\beta_{jk} & =\frac{(q^{\alpha+1},-c,-c^{-1}q;q)_{\infty}}{(q,-cq^{\alpha+1},-c^{-1}q^{-\alpha};q)_{\infty}}\frac{q^{\alpha(j+k)+j^{2}+k^{2}}}{(q;q)_{j}(q;q)_{k}}\label{eq:3.32}\\
\times & \sum_{m=\max(j,k)}^{n}\frac{(q;q)_{m}q^{m}}{(q^{\alpha+1};q)_{m}}\frac{(q^{\alpha+j+1};q)_{m-j}(q^{\alpha+k+1};q)_{m-k}}{(q;q)_{m-j}(q;q)_{m-k}},\nonumber \end{align}
 then, the inverse matrix $(\gamma_{jk})_{j,k=0}^{n}$ of \begin{align}
 & \left((q^{\alpha+1};q)_{j+k}q^{-\binom{j+k}{2}}\right)_{j,k=0}^{n}\label{eq:3.33}\end{align}
 has element\begin{align}
\gamma_{jk} & =\frac{q^{j^{2}-j+k^{2}-k}}{(q;q)_{j}(q;q)_{k}}\sum_{m=\max(j,k)}^{n}\frac{(q;q)_{m}q^{m}}{(q^{\alpha+1};q)_{m}}\frac{(q^{\alpha+1};q)_{m-j}(q^{\alpha+1};q)_{m-k}}{(q;q)_{m-j}(q;q)_{m-k}}\label{eq:3.34}\\
= & \frac{q^{[j(j-1)+k(k-1)]}}{(q^{\alpha+1};q)_{j}(q^{\alpha+1};q)_{k}}\sum_{m=\max(j,k)}^{n}\frac{(q^{\alpha+1};q)_{m}q^{m}}{(q;q)_{m}}\left[\begin{array}{c}
m\\
j\end{array}\right]_{q}\left[\begin{array}{c}
m\\
k\end{array}\right]_{q},\nonumber \end{align}
 or\begin{align}
\gamma_{jk} & =\frac{q^{[j(j-1)+k(k-1)]}}{(q^{\alpha+1};q)_{j}(q^{\alpha+1};q)_{k}}\sum_{m=\max(j,k)}^{n}\frac{(q^{\alpha+1};q)_{m}q^{m}}{(q;q)_{m}}\left[\begin{array}{c}
m\\
j\end{array}\right]_{q}\left[\begin{array}{c}
m\\
k\end{array}\right]_{q}.\label{eq:3.35}\end{align}
 Since the formulas \eqref{eq:3.31}, \eqref{eq:3.33} and \eqref{eq:3.35}
only involve rational functions of $q^{\alpha}$ and $q^{k},k\in\mathbb{N}$,
thus the restriction $q$ may be dropped to get Theorem \ref{thm:2.9}.

\subsection{Proof for Theorem \ref{thm:2.10}}

The discrete $q$-Hermite polynomials II are defined as \cite{Koekoek}\begin{align}
\tilde{h}_{n}(x;q): & =x^{n}{}_{2}\phi_{n}\left(\begin{array}{c}
q^{-n},q^{-n+1}\\
0\end{array};q^{2};-\frac{q^{2}}{x^{2}}\right)\label{eq:81}\end{align}
 for $n\ge0$ and we assume that\begin{align}
\tilde{h}_{-1}(x;q): & =0,\quad\tilde{h}_{0}(x;q):=1.\label{eq:3.37}\end{align}
 The polynomials $\left\{ \tilde{h}_{n}(x;q)\right\} _{n=0}^{\infty}$
satisfy the difference equation\begin{align}
D_{q^{-1}}^{j}\tilde{h}_{n}(x;q) & =\frac{(q;q)_{n}q^{-jn+j(j+1)/2}}{(1-q)^{j}(q;q)_{n-j}}\tilde{h}_{n-j}(x;q).\label{eq:3.38}\end{align}
 They also satisfy the orthogonality\begin{align}
\sum_{k=-\infty}^{\infty}\left[\tilde{h}_{m}(cq^{k};q)\tilde{h}_{n}(cq^{k};q)+\tilde{h}_{m}(-cq^{k};q)\tilde{h}_{n}(-cq^{k};q)\right]w(cq^{k})q^{k} & =h_{n}(c|q)\delta_{mn}\label{eq:3.39}\end{align}
 for $m,n\ge0$, where\begin{align}
w(x): & =\frac{1}{(-x^{2};q^{2})_{\infty}},\label{eq:3.40}\end{align}
 and\begin{align}
h_{n}(c|q): & =2\frac{(q^{2},-c^{2}q,-c^{-2}q;q^{2})_{\infty}}{(q,-c^{2},-c^{-2}q^{2};q^{2})_{\infty}}\frac{(q;q)_{n}}{q^{n^{2}}}\label{eq:3.41}\end{align}
 for some $c>0$. thus the orthonormal polynomial\begin{align}
h_{n}(x): & =\frac{\tilde{h}_{n}(x;q)}{\sqrt{h_{n}(c|q)}}\label{eq:3.42}\end{align}
has the leading coefficient\begin{align}
\gamma_{n} & =\frac{1}{\sqrt{h_{n}(c|q)}}.\label{eq:3.43}\end{align}
 The moments of the measure could be calculated via formula \ref{eq:2.7},\begin{align}
\mu_{n} & =\sum_{k=-\infty}^{\infty}\left[(cq^{k})^{n}+(-cq^{k})^{n}\right]w(cq^{k})q^{k}\label{eq:3.44}\\
= & \left[1+(-1)^{n}\right]c^{n}\sum_{k=-\infty}^{\infty}\frac{q^{(n+1)k}}{(-c^{2}q^{2k};q^{2})_{\infty}}\nonumber \\
= & \frac{\left[1+(-1)^{n}\right]c^{n}}{(-c^{2};q^{2})_{\infty}}\sum_{k=-\infty}^{\infty}(-c^{2};q^{2})_{k}q^{(n+1)k}\nonumber \\
= & \frac{\left[1+(-1)^{n}\right]c^{n}}{(-c^{2};q^{2})_{\infty}}\frac{(q^{2},-c^{-2}q^{-n+1},-c^{2}q^{n+1};q^{2})}{(-c^{-2}q^{2},q^{n+1};q^{2})_{\infty}}\nonumber \\
= & \frac{\left[1+(-1)^{n}\right]q^{-n^{2}/4}(q^{2},-c^{-2}q,-c^{2}q;q^{2})_{\infty}}{(-c^{2},-c^{-2}q^{2},q^{n+1};q^{2})_{\infty}}.\nonumber \end{align}
 Thus, for $n=0,1,...$, the matrix\begin{align}
H_{n} & =\left(\frac{\left[1+(-1)^{j+k}\right]q^{-(j+k)^{2}/4}(q^{2},-c^{-2}q,-c^{2}q;q^{2})_{\infty}}{(-c^{2},-c^{-2}q^{2},q^{j+k+1};q^{2})_{\infty}}\right)_{j,k=0}^{n}\label{eq:3.45}\end{align}
 has determinant\begin{align}
\det H_{n} & =\frac{\prod_{k=0}^{n}(q;q)_{k}}{q^{n(n+1)(2n+1)/6}}\left[2\frac{(q^{2},-c^{2}q,-c^{-2}q;q^{2})_{\infty}}{(q,-c^{2},-c^{-2}q^{2};q^{2})_{\infty}}\right]^{n+1}.\label{eq:3.46}\end{align}
 Therefore,\begin{align}
\det\left(\frac{\left[1+(-1)^{j+k}\right]}{2}q^{-(j+k)^{2}/4}(q;q^{2})_{\left\lfloor \frac{j+k}{2}\right\rfloor }\right)_{j,k=0}^{n} & =\frac{\prod_{k=0}^{n}(q;q)_{k}}{q^{n(n+1)(2n+1)/6}}.\label{eq:3.47}\end{align}
 For $n=0,1,...$, the matrix $H_{n}^{-1}=(\beta_{jk})_{j,k=0}^{n}$
has element\begin{align}
\beta_{jk} & =(-1)^{j+k}q^{j^{2}+k^{2}}\frac{(q,-c^{2},-c^{-2}q^{2};q^{2})_{\infty}}{2(q^{2},-c^{2}q,-c^{-2}q;q^{2})_{\infty}}\label{eq:3.48}\\
\times & \sum_{m=\max(j,k)}^{n}\frac{q^{m^{2}-(j+k)m}\tilde{h}_{m-j}(0;q)\tilde{h}_{m-k}(0;q)}{(q;q)_{m}}\left[\begin{array}{c}
m\\
j\end{array}\right]_{q}\left[\begin{array}{c}
m\\
k\end{array}\right]_{q}.\nonumber \end{align}
From the generating function\begin{align}
\frac{(-xt;q)_{\infty}}{(-t^{2};q^{2})_{\infty}} & =\sum_{n=0}^{\infty}\frac{q^{\binom{n}{2}}t^{n}}{(q;q)_{n}}\tilde{h}_{n}(x;q),\label{eq:3.49}\end{align}
we get\begin{align}
\tilde{h}_{n}(0;q) & =\frac{(q;q^{2})_{\infty}\cos\frac{n\pi}{2}}{(q^{1+n};q^{2})_{\infty}q^{\binom{n}{2}}}.\label{eq:3.50}\end{align}
Then, the matrix\begin{align}
 & \left(\frac{\left[1+(-1)^{j+k}\right]}{2}q^{-(j+k)^{2}/4}(q;q^{2})_{\left\lfloor \frac{j+k}{2}\right\rfloor }\right)_{j,k=0}^{n}\label{eq:3.51}\end{align}
 has inverse matrix\begin{align}
 & \left(\sum_{m=\max(j,k)}^{n}\frac{\left[\begin{array}{c}
m\\
j\end{array}\right]_{q}\left[\begin{array}{c}
m\\
k\end{array}\right]_{q}(q,q^{2})_{\left\lfloor \frac{m-j}{2}\right\rfloor }(q,q^{2})_{\left\lfloor \frac{m-k}{2}\right\rfloor }\cos\frac{(m-j)\pi}{2}\cos\frac{(m-k)\pi}{2}}{(-1)^{m}(q;q)_{m}q^{-m-\binom{j}{2}-\binom{k}{2}}}\right)_{j,k=0}^{n}.\label{eq:3.52}\end{align}
Since formulas \eqref{eq:3.47}, \eqref{eq:3.51}and \eqref{eq:3.52}
contain only rational functions of $q^{k},k\in\mathbb{N}$, thus the
restriction may be dropped and Theorem \ref{thm:2.10} follows.

\begin{acknowledgement*}
This work is partially supported by Chinese National Natural Science
Foundation grant No.10761002, Guangxi Natural Science Foundation grant
No.0728090.
\end{acknowledgement*}

\end{document}